\documentclass[12pt,a4paper]{article}
\usepackage{amssymb,amsmath}
\usepackage{amsfonts}
\usepackage{a4}
\newcommand{\para}{\par\vspace{.25cm}}

\newtheorem{theorem}{Theorem}
\newtheorem{lemma}{Lemma}

\begin{document}
\baselineskip 18pt \title{\bf The algebraic structure of finite metabelian group algebras}
\author{ Gurmeet K. Bakshi\\ {\em \small Centre for Advanced Study in
Mathematics}\\
{\em \small Panjab University, Chandigarh 160014, India}\\{\em
\small email: gkbakshi@pu.ac.in} \and Shalini Gupta\\
{\em \small Department of Mathematics}\\ {\em \small Punjabi University, Patiala 147002, India}\\
{\em \small email: {gupta\_math@yahoo.com}}
\and Inder Bir S. Passi \\  {\em \small Centre for Advanced Study
in Mathematics}\\ {\em \small Panjab University, Chandigarh 160014, India} \\
{\small \& }\\
{\em \small Indian Institute of Science Education and Research}\\
{\em \small Sector 81, Mohali 140306, India}\\{\em \small email: ibspassi@yahoo.co.in } }
\date{} \maketitle
 \date{}
 {\maketitle}
\begin{abstract}\noindent
{ An algorithm for the explicit computation of a complete set of primitive
 central idempotents, Wedderburn decomposition and the automorphism group of the  semisimple group algebra of a 
 finite metabelian group is developed. The algorithm is illustrated with its application to the
semisimple group algebra of an arbitrary metacyclic group, and  certain indecomposable groups whose central
quotient is the Klein four-group.}
\end{abstract}\vspace{.25cm}

  \vspace{.25cm} 
\section{Introduction} Let  $\mathbb{F}_{q}[G]$ be the group algebra of a finite group $G$ over the finite field
$\mathbb{F}_{q}$ with $q$ elements.  In order to understand the algebraic structure of  $\mathbb{F}_{q}[G]$,
in the semisimple case, i.e., when $q$ is coprime to the order of $G$, an essential step is to compute a complete set 
of primitive central idempotents and the Wedderburn
decomposition of  $\mathbb{F}_{q}[G]$. These computations, in turn, help to investigate the automorphism group and the 
unit group of  $\mathbb{F}_{q}[G]$.
\para Let $(H,\, K)$ be a {\it strongly Shoda pair} (\cite{bro1}, Definition 5) of $G$ and let $C$ be a $q$-cyclotomic 
coset of  $\operatorname{Irr}(K/H),$ the set of irreducible characters of $K/H$ over the algebraic closure 
$\overline{\mathbb{F}}_{q}$ of $\mathbb{F}_{q},$  
corresponding to a generator of $\operatorname{Irr}(K/H)$.
Broche and Rio (\cite{bro1}, Theorem 7 )  proved that the pair $((H,\, K),\, C)$ defines a primitive central
idempotent,
 $e_{C}(G,\, K,\, H)$, of $\mathbb{F}_{q}[G]$. They further proved that if $G$ is an
 abelian-by-supersolvable group, then every primitive central idempotent
of the semisimple group algebra $\mathbb{F}_{q}[G]$ is defined by  a pair of the type $((H,\, K), \,C)$.
 However, it is possible that two distinct such pairs may define the same primitive central idempotent. Thus
in order to determine the algebraic structure of  $\mathbb{F}_{q}[G],\,G$ abelian-by-supersolvable, the problem  lies 
in finding a set
$\mathfrak{S}$ of  pairs of type $((H,\, K),\, C)$ so that  $\{e_{C}(G,\,K,\,H) \mid ((H,\,K),\, C) \in
\mathfrak{S} \}$ is a complete irredundant set of primitive central idempotents of $\mathbb{F}_{q}[G]$.
\para In Section 2, we provide an efficient  algorithm for the solution of this problem for $\mathbb{F}_{q}[G],$
when $G$ is metabelian (Theorem 2). Our analysis, in turn,
leads to an explicit description of the Wedderburn
decomposition and the automorphism group of  $\mathbb{F}_{q}[G]$ (Theorem 3).
   In Section 3, we illustrate our algorithm with its application to an arbitrary finite metacyclic group $G = \langle 
   a, b \,|\, a^{n} = 1,\, b^{t} = a^{k},
\, b^{-1}ab=
a^{r}\rangle,$ $n,\,t,\,k,\,r$  natural numbers, $r^{t}\equiv 1~({\rm mod}\, n),\,k(r-1)\equiv 0~({\rm mod}\,n$),
 of order $nt$ coprime to $q,$ and provide an alternative way of finding a complete set of the primitive central 
 idempotents to those given in
\cite{sha2}. We next apply our result, in Section 4, to indecomposable groups  $G$ whose central quotient, $G/Z(G),$ is 
the Klein four-group.
  It is known, (\cite{goo1}, Chapter 5), that such groups breakup into five different classes.
 Ferraz, Goodaire and Milies \cite{fer2} have given, in each case, a lower bound on the number of simple components of
the semisimple finite group algebra $\mathbb{F}_{q}[G].$ We provide the complete algebraic structure of the semisimple 
group
algebra $\mathbb{F}_{q}[G]$ for group $G$ in two of the five classes; thus improving  Theorems 3.1 \& 3.2 of
 \cite{fer2}.
\section{Metabelian groups}
Let $G$ be a finite group. We adopt the following notation:
$$\begin{array}{ll} 
  e_{\mathbb{F}_{q}}(\chi) & \frac{\chi(1)}{|G|} \sum_{ \atop{\sigma \in
\operatorname{Gal}(\mathbb{F}_{q}(\chi)/\mathbb{F}_{q})}}\sigma(\chi(g))g^{-1},~
\mbox{the primitive central idempotent of}~\mathbb{F}_{q}[G]\\&~\mbox{ determined by }~
\chi\in \operatorname{Irr}(G), \mbox{where}~\mathbb{F}_{q}(\chi)~\mbox{ is the field obtained by adjoining to}\\&
\mathbb{F}_{q},~\mbox{all the character values}~\chi(g), g \in G ;\vspace{2mm}\\
\mathfrak{A}& \mbox{the set of pairs}~(H,\,K),~\mbox{where}~H\unlhd K\leq G~\mbox{and}~ K/H~\mbox{is cyclic};
\vspace{2mm}\\
 \mathcal{C}(K/H)& \mbox{the set of}~q\mbox{-cyclotomic cosets of}~\operatorname{Irr}(K/H)~\mbox{containing
the generators of}~\\&\operatorname{Irr}(K/H),~\mbox{where}~(H,\,K)\in \mathfrak{A};
\vspace{2mm}\\\mathcal{R}(K/H)& \mbox{the set of distinct orbits of}~\mathcal{C}(K/H)~\mbox{under the action of}~
N_{G}(H)\cap N_{G}(K)\\&
~\mbox{on}~\mathcal{C}(K/H)\mbox{ given by}~g.C=C^{g}(:=g^{-1}Cg),~g\in N_{G}(H)\cap N_{G}(K),\\& C\in \mathcal{C}(K/H),
\mbox{where}~N_{G}(H)~\mbox{denotes the normalizer of}~H\mbox~{in}~G;\vspace{2mm}\\
E_{G}(K/H)&\mbox{the stabilizer of any}~C\in \mathcal{C}(K/H)~\mbox{under the above action of}\\&
N_{G}(H)\cap N_{G}(K)
~\mbox{on}~\mathcal{C}(K/H)~(\mbox{note that this stabilizer does not depend on}~C);\vspace{2mm}\\
\varepsilon_{C}(K, H) &|K|^{-1}\sum_{g \in K}
\operatorname{tr}_{\mathbb{F}_{q}(\zeta)/\mathbb{F}_{q}}(\chi(\overline{g}))g^{-1},~\mbox{where}~\chi ~\mbox{is a 
representative of the}\\&  q\mbox{-cyclotomic coset}~C,~\mbox{and}~\zeta~\mbox{is a primitive}~[K : H]\mbox{-th root of 
unity in}~ \overline{\mathbb{F}}_{q},\\ &(H,\,K)\in \mathfrak{A},\,C\in \mathcal{C}(K/H); \vspace{2mm}\\
e_{C}(G,\,K,\,H)& \mbox{the sum of distinct}~G\mbox{-conjugates of}~\varepsilon_{C}(K, H).
\end{array}$$\par
For the rest of this section, we assume $G$ to be a finite metabelian group.

\subsection{Primitive central idempotents}  We follow
 the notation introduced in \cite{bkp}.
 Let
\begin{itemize}\item $\mathcal{A}$:= a fixed maximal abelian subgroup of $G$ containing its derived subgroup $G'.$
 \item $\mathcal{T}$ := the set of all subgroups $D$ of $G$ with $D\leq \mathcal{A}$ and $\mathcal{A}/D$  cyclic.
\end{itemize}
For $D_{1},\,D_{2}\in
\mathcal{T},$ we say that $D_{1}$ is equivalent to $D_{2}$ if there exists $g\in G$ such that $D_{2}= D_{1}^{g}.$ 
Let
\begin{itemize}
 \item $\mathcal{T}_{G}$:= a set of representatives
of the distinct equivalence classes of $\mathcal{T}.$
\end{itemize}
  For $D\in \mathcal{T},$ let
\begin{itemize}
 \item $K_{D}$ := a fixed maximal element of
$\{K~|~\mathcal{A} \leq K \leq G,~K'\leq D\}.$
\item $\mathcal{R}(D)$ := the set of those
linear representations of $K_{D}$ over $\overline{\mathbb{F}}_{q}$ whose restriction to $\mathcal{A}$ has kernel 
$D.$
\item $\mathcal{R}_{C}(D)$ := a complete set of those representations in $\mathcal{R}(D)$ which are
not mutually $G$-conjugate.
\end{itemize}
\par\vspace{.35cm}
The following result is proved in \cite{bkp} for complex irreducible representations. However, the analogous
proof works for the irreducible representations of
$G$ over $\overline{\mathbb{F}}_{q}.$ \vspace{.15cm}
\begin{theorem}{\rm \cite{bkp}}\label{p2} Let $G$ be a finite metabelian group with $\mathcal{A}$  and 
$\mathcal{T}_{G}$ as defined
above. Then
$$\Omega=\{ \rho^{G}\,|\, \rho \in \mathcal{R}_{C}(D),~D\in \mathcal{T}_{G}\},$$ is a complete  set of  inequivalent 
irreducible representations of $G$ over $\overline{\mathbb{F}}_{q},$ where $\rho^{G}$ denotes $\rho$ induced to 
$G$.\par
 Furthermore, $\rho^{G}\in \Omega$ is faithful if, and only if, $D$ is core-free in $G$, i.e., $\cap_{x \in G}D^{x} =
  \{1\}.$\end{theorem}
\par\vspace{.25cm}
  For  $N \unlhd G$ with $$\mathcal{A}_{N}/N :=~\mbox{ a maximal abelian subgroup of}~ G/N~\mbox{ 
  containing}~{(G/N)}',$$  define
$$\mathcal{S}_{G/N} := \{(D/N,\,\mathcal{A}_{N}/N)\,|\,D/N\in \mathcal{T}_{G/N},\,
D/N ~ \mbox{core-free in}~G/N\}.$$
Let \[\mathcal{S}:= \{(N,\,D/N,\,\mathcal{A}_{N}/N)\,|\,N \unlhd G,\,\mathcal{S}_{G/N}\neq
\o{},\,
(D/N,\,\mathcal{A}_{N}/N)\in \mathcal{S}_{G/N}\}.\]\\ \vspace{.15cm}
 By (\cite{bkp}, Lemma 6), each element $(N,\,D/N,\,\mathcal{A}_{N}/N)\in
\mathcal{S}$ defines a strongly Shoda pair $(D,\,\mathcal{A}_{N})$ in $G$ and the mapping 
$((N,\,D/N,\,\mathcal{A}_{N}/N) \mapsto (D,\,\mathcal{A}_{N})$ from $\mathcal{S}$ to the set of strongly Shoda pairs of 
$G$ is one-one. Thus $\mathcal{S}$ may be regarded as a subset of strongly Shoda pairs of $G.$\par
% \[\mathfrak{S}:= \{((N,\,D/N,\,\mathcal{A}_{N}/N),\,C)\,|\,
% (N,\,D/N,\,\mathcal{A}_{N}/N)\in \mathcal{S},\,C\in R(\mathcal{A}_{N}/D)\}.\]
Let $\rm{ord}_{n}(q)$ denote the order of $q$ modulo $n,\,n\geq 1.$ We prove the following:
\begin{theorem}\label{T1}Let  $\mathbb{F}_{q}$ be a finite field with $q$ elements and
 $G$ a finite metabelian group. Suppose that $\gcd(q,\,|G|) = 1.$ Then,
\begin{quote}
 $(i)$~ $\{e_{C}(G,\,\mathcal{A}_{N},\,D)\,\mid\,(N,\,D/N,\,\mathcal{A}_{N}/N)\in \mathcal{S},\,C\in 
 \mathcal{R}(\mathcal{A}_{N}/D)\}$ is a complete
set of primitive central idempotents of $\mathbb{F}_{q}[G];$ \par\vspace{.25cm}

$(ii)$~for  $(N,\,D/N,\,\mathcal{A}_{N}/N)\in \mathcal{S}$ and $C \in \mathcal{R}(\mathcal{A}_{N}/D),$ the simple 
component $\mathbb{F}_{q}[G] e_{C}(G,\,\mathcal{A}_{N},\,D)$ is isomorphic to 
$M_{[G:\mathcal{A}_{N}]}(\mathbb{F}_{q^{o(\mathcal{A}_{N},\,D) }}),$ the algebra of $[G:\mathcal{A}_{N}] \times 
[G:\mathcal{A}_{N}]$ matrices over the field $\mathbb{F}_{q^{o(\mathcal{A}_{N},\,D) }}$, 
where $ o(\mathcal{A}_{N},\,D)=\frac{{\rm 
ord}_{([\mathcal{A}_{N}:D])}(q)}{[E_{G}(\mathcal{A}_{N}/D):\mathcal{A}_{N}]}.$
 \end{quote}

\end{theorem}
{\bf{Proof}.} $(i)$
~Let \begin{equation}\label{b1}\mathfrak{S}:= \{((N,\,D/N,\,\mathcal{A}_{N}/N),\,C)\,|\,
(N,\,D/N,\,\mathcal{A}_{N}/N)\in \mathcal{S},\,C\in R(\mathcal{A}_{N}/D)\}.\end{equation}
 If $((N,\,D/N,\,\mathcal{A}_{N}/N),\,C)
\in \mathfrak{S},$ then, by (\cite{bkp}, Lemma 6), $(D,\,\mathcal{A}_{N})$ is a strongly Shoda pair in $G,$ and 
therefore, by
 (\cite{bro1}, Theorem 7 ), $e_{C}(G,\,\mathcal{A}_{N},\,D)$ is a
primitive central idempotent of $\mathbb{F}_{q}[G].$ Thus we have a map
$$\pi: ((N,\,D/N,\,\mathcal{A}_{N}/N),\,C)\mapsto e_{C}(G,\,\mathcal{A}_{N},\,D)$$ from  $\mathfrak{S}$ to a complete
 set of primitive central idempotents of $\mathbb{F}
_{q}[G].$
In order to prove the Theorem, we need
to prove that $\pi$ is 1-1 and onto.\par\vspace{.25cm}
To show that $\pi$ is onto, let $e$ be a primitive central idempotent of $\mathbb{F}_{q}[G].$ We have $e= 
e_{\mathbb{F}_{q}}(\chi),$ for some
$\chi\in \operatorname{Irr}(G).$ Let $\tau$ be a representation affording $\chi$ and $N =
\operatorname{ker}\tau,$ the kernel of the character $\tau$.
Let $\overline{\tau}$ be the corresponding faithful representation of $G/N.$ By Theorem \ref{p2}, it follows that
there exists a unique pair $(D/N,\,\mathcal{A}_{N}/N)\in S_{G/N}$ and a representation
$\overline{\rho}$ of
$\mathcal{A}_{N}/N$ with kernel $D/N$ such that $\overline{\tau} = {\overline{\rho}}^{G/N}.$ This yields $\chi =
 \psi^{G},$ where $\psi$ is the character afforded by $\rho : \mathcal{A}_{N}\rightarrow \overline{\mathbb{F}}_{q}$ 
 given by
$\rho(x) = \overline{\rho}(xN)$. Since $\operatorname{ker}\psi = D$, by (\cite{sha2}, Lemma 1), we have
\begin{equation}\label{E4} e_{\mathbb{F}_{q}}(\chi) = e_{C}(G,\,\mathcal{A}_{N},\,D),\end{equation}
 where $C\in \mathcal{R}(\mathcal{A}_{N}/D)$ is
the $q$-cyclotomic coset of
$\overline{\psi}$ and consequently $\pi$ is onto.\par\vspace{.25cm}
To show that $\pi$ is 1-1, let $((N,\,D/N,\,\mathcal{A}_{N}/N),\,C)
$ and\linebreak $((\tilde{N},\,\tilde{D}/\tilde{N},\,\mathcal{A}_{\tilde{N}}/\tilde{N}),\,\tilde{C}) \in \mathfrak{S} $ 
be such that

\begin{equation}\label{e4}e_{C}(G,\,\mathcal{A}_{N},\,D) =
e_{\tilde{C}}(G,\,\mathcal{A}_{\tilde{N}},\,\tilde{D}).\end{equation} Let $\rho \in \mathcal{R}_{C}(D)$, $
\tilde{\rho} \in \mathcal{R}_{\tilde{C}}(\tilde{D})$ and
$\chi$ and $\tilde{\chi}$  be the character afforded by  $\rho^{G}$ and $\tilde{\rho}^{G}$ respectively.
 By (\cite{sha2}, Lemma 1),
 $ e_{\mathbb{F}_{q}}(\chi) = e_{C}(G,\,\mathcal{A}_{N},\,D)
$ and $ e_{\mathbb{F}_{q}}(\tilde{\chi}) = e_{\tilde{C}}(G,\,\mathcal{A}_{\tilde{N}},\,\tilde{D}).$
Therefore, equation (\ref{e4}) implies that  $e_{\mathbb{F}_{q}}(\chi)=
 e_{\mathbb{F}_{q}}(\tilde{\chi}),$ which, in turn, implies that
\begin{equation}\label{E3}
 \tilde{\chi} = \sigma \circ \chi,~\sigma\in \operatorname{Gal}(\mathbb{F}_{q}(\chi)/\mathbb{F}_{q}).
\end{equation} Consequently, $\tilde{N} = \operatorname{ker}(\tilde{\chi})=  \operatorname{ker}(\chi)= N.$
  Also, by going modulo $N,$ it follows from equation (\ref{E3}) and Theorem \ref{p2},  that $D/N$ and $\tilde{D}/N$ 
  are conjugate in $G/N.$ This gives
$D/N = \tilde{D}/N,$ i.e., $D = \tilde{D} $.
Next, if  $\{z_{1},\,z_{2},\,\ldots,\,z_{k}\}$ is a transversal of
$E_{G}(\mathcal{A}_{N}/D)$ in $G,$ then,  by (\cite{bro1}, Lemma 4 ) and equation (\ref{e4}), we have
\begin{equation}\label{e3}\displaystyle\sum_{j=1}^{k}\varepsilon_{C^{z_{j}}}(\mathcal{A}_{N},\,D^{z_{j}})=
\displaystyle\sum_{j=1}^{k}\varepsilon_{\tilde{C}^{z_{j}}}(\mathcal{A}_{N},\,D^{z_{j}}).
\end{equation}
Since both the sides of the above equation are primitive central idempotents in $\mathbb{F}_{q}[\mathcal{A}_{N}],$ it 
follows that,
for some $j,~1\leq j\leq k,$
\begin{equation}\label{e5}
 \varepsilon_{C}(\mathcal{A}_{N},\,D)=\varepsilon_{\tilde{C}^{z_{j}}}
(\mathcal{A}_{N},\,D^{z_{j}}).\end{equation}
However, by (\cite{bro1}, Proposition 2), $\varepsilon_{C}(\mathcal{A}_{N},\,D) = e_{\mathbb{F}_{q}}(\rho),$
 and
$\varepsilon_{\tilde{C}^{z_{j}}}
(\mathcal{A}_{N},\,\tilde{D}^{z_{j}})=  e_{\mathbb{F}_{q}}(\tilde{\rho}^{z_{j}}).$ Therefore, we have by equation 
(\ref{e5}),
$e_{\mathbb{F}_{q}}(\rho)= e_{\mathbb{F}_{q}}(\tilde{\rho}^{z_{j}}),$ which, as before, gives
$ D = \operatorname{ker}\rho =
\operatorname{ker}\tilde{\rho}^{z_{j}} = \tilde{D}^{z_{j}}= D^{z_{j}},$ i.e., $z_{j}\in N_{G}(D).$
Consequently, $C$ and $C'$ have same orbits. This proves that $\pi$ is
1-1.  $\Box$ \vspace{.3cm} \\ $(ii)$ follows from (\cite{bro1}, Corollary $9$). $\Box$
\subsection{Wedderburn decomposition and automorphism group} We continue with the notation introduced in $\S 2.1.$ Let 
$\operatorname{Aut}(\mathbb{F}_{q}[G])$ be the group of $\mathbb{F}_{q}$-automorphisms of $\mathbb{F}_{q}[G]$. For $ n \geq 1$, let 
$\mathbb{Z}_{n}$ be the additive group of integers modulo $n$,  $S_{n}$  the symmetric group of degree $n$, $SL_{n} 
(K)$ the group of invertible $ n \times n$ matrices over the field $K$ of determinant $1,$ and for any algebra $K,$ let 
 $K^{(n)}$ denote the direct sum of $n$ copies of $K.$ Let $\xi$ be a primitive $|G|$-th root
of unity in  $\overline{\mathbb{F}}_{q}$. Let $(N,\,D/N,\,\mathcal{A}_{N}/N)\in \mathcal{S}.$ Then
$\mathcal{A}_{N}/D$ is a cyclic
group generated by $aD,$ say. Let $x_{1},\,x_{2},\ldots,\,x_{t}$ be a transversal of  $\mathcal{A}_{N}$ in $G,$ and 
$r_{i},\,
1\leq i\leq t,$ be integers such that $x_{i}^{-1}ax_{i}D = a^{r_{i}}D.$ Let
 $\zeta ={\xi}^{|G|/[\mathcal{A}_{N}:D]}$, and 
 $\mathcal{K}(N,\,D/N,\,\mathcal{A}_{N}/N)$\index{$\mathcal{K}(N,\,D/N,\,\mathcal{A}_{N}/N)$} be the subfield of
$\overline{\mathbb{F}}_{q}$ obtained by adjoining the $t$ elements $\displaystyle\sum_{i=1}^{t}\zeta^{jr_{i}},~
1\leq j\leq t-1$ to $\mathbb{F}_{q}.$ It is easily seen that the field
$ \mathcal{K}(N,\,D/N,\,\mathcal{A}_{N}/N)$ is independent of the choice
of transversal of $\mathcal{A}_{N}$ in $G.$ \par\vspace{.25cm}
For $d|[G:G']$ and $l|[\mathbb{F}_{q}(\xi):\mathbb{F}_{q}],$ let
 $\mathcal{S}_{d,\,l}$\index{$\mathcal{S}_{d,\,l}$} be the set of those $(N,\,D/N,\,\mathcal{A}_{N}/N)\in \mathcal{S}$ 
 such that
\begin{quote} $(i)$~$[G:\mathcal{A}_{N}]=d,$ \vspace{2mm}\\
$(ii)$~$[\mathcal{K}(N,\,D/N,\,\mathcal{A}_{N}/N):\mathbb{F}_{q}] = l.$
 \end{quote} Clearly
$\mathcal{S}_{d,\,l},~d|[G:G'],~l|[\mathbb{F}_{q}(\xi):\mathbb{F}_{q}],$
 are disjoint and  $\mathcal{S} = \bigcup\limits_{d\,|\,[G:G']\atop{ l\,|\,[\mathbb{F}_{q}(\xi):
\mathbb{F}_{q}]}}\mathcal{S}_{d,\,l}.$\vspace{.25cm}
 \begin{theorem}\label{T2} With the above notation,\\ \\
$(i)$~  $\mathbb{F}_{q}[G]\cong \bigoplus\limits_{d \mid [G:G']\atop
{ l \mid [\mathbb{F}_{q}(\xi):
\mathbb{F}_{q}]}} {M_{d}(\mathbb{F}_{q^{l}})}^{(\alpha_{d,\,l})},$ \\ \\
$(ii)$~ $\mathrm{Aut}(\mathbb{F}_{q}[G])\cong \bigoplus\limits_{d\,|\,[G:G']\atop
{ l\,|\,[\mathbb{F}_{q}(\xi):
\mathbb{F}_{q}]}}{K}^
{(\alpha_{d,\,l})}_{d,\,l}\rtimes S_{\alpha_{d,\,l}},$\\ \\  where $K _{d,\,l}=
\mathrm{SL}_{d}(\mathbb{F}_{q^{l}})
\rtimes \mathbb{Z}_{l},$ a semi direct product of $\mathrm{SL}_{d}(\mathbb{F}_{q^{l}})$ by 
$\mathbb{Z}_{l}$,  and $\alpha_{d,\,l} =
 \displaystyle\sum_{(N,\,D/N,\,\mathcal{A}_{N}/N) \in S_{d,\,l}}| \mathcal{R}(\mathcal{A}_{N}/D)|.$

\end{theorem}\vspace{.25cm}
{\bf{Proof.}}~ $(i)$ It follows from Theorem \ref{T1} that for
$((N,\,D/N,\,\mathcal{A}_{N}/N),C)$ $\in \mathfrak{S},$ where $\mathfrak{S}$ is as defined in equation (\ref{b1}),
\[\mathbb{F}_{q}[G] e_{C}(G,\,\mathcal{A}_{N},\,D)\cong 
M_{[G:\mathcal{A}_{N}]}(\mathbb{F}_{q^{o(\mathcal{A}_{N},\,D)}})
\]  Thus we have,
\begin{eqnarray*}\mathbb{F}_{q}[G] &\cong& \bigoplus_{((N,\,D/N,\,\mathcal{A}_{N}/N),\,C)\in
\mathfrak{S}}
 \mathbb{F}_{q}[G] e_{C}(G,\,\mathcal{A}_{N},\,D)\\ \\ &\cong& \bigoplus_{(N,\,D/N,\,\mathcal{A}_{N}/N)\in
\mathcal{S}}\bigoplus_{C\in \mathcal{R}(\mathcal{A}_{N}/D)}
 \mathbb{F}_{q}[G] e_{C}(G,\,\mathcal{A}_{N},\,D)\\ \\ &\cong&
\bigoplus_{(N,\,D/N,\,\mathcal{A}_{N}/N)\in
\mathcal{S}}\bigoplus_{C\in \mathcal{R}(\mathcal{A}_{N}/D)} M_{[G:\mathcal{A}_{N}]}
(\mathbb{F}_{q^{o(\mathcal{A}_{N},\,D)}})\\ \\ &\cong& \bigoplus\limits_{d \mid [G:G']\atop
{ l \mid [\mathbb{F}_{q}(\xi):
\mathbb{F}_{q}]}} \bigoplus_{(N,\,D/N,\,\mathcal{A}_{N}/N)\in
\mathcal{S}_{d,\,l}}\bigoplus_{C\in \mathcal{R}(\mathcal{A}_{N}/D)} M_{[G:\mathcal{A}_{N}]}
(\mathbb{F}_{q^{o(\mathcal{A}_{N},\,D)}})\\ \\ &\cong& \bigoplus\limits_{d \mid [G:G']\atop
{ l \mid [\mathbb{F}_{q}(\xi):
\mathbb{F}_{q}]}} \bigoplus_{(N,\,D/N,\,\mathcal{A}_{N}/N)\in
\mathcal{S}_{d,\,l}}{ M_{[G:\mathcal{A}_{N}]}
(\mathbb{F}_{q^{o(\mathcal{A}_{N},\,D)}})}^{(|\mathcal{R}(\mathcal{A}_{N}/D)|)}\end{eqnarray*} For $d\mid 
[G:G'],~l\mid
[\mathbb{F}_{q}(\xi):\mathbb{F}_{q}],$ and $(N,\,D/N,\,\mathcal{A}_{N}/N)\in
\mathcal{S}_{d,\,l},$ we show that \begin{equation}\label{E5}o(\mathcal{A}_{N},\,D) =
[\mathcal{K}(N,\,D/N,\,\mathcal{A}_{N}/N):\mathbb{F}_{q}] = l.\end{equation}
\par\vspace{.25cm}  If $\rho \in \mathcal{R}_{C}(D)$ and $\chi$ is the character afforded by $\rho^{G},$ then, by
(\cite{sha2}, Lemma 1),
$$[E_{G}(\mathcal{A}_{N}/D):\mathcal{A}_{N}] = [\mathbb{F}_{q}(\zeta):\mathbb{F}_{q}(\chi)].$$  However, note that
 $$\mathbb{F}_{q}(\chi) =
\mathcal{K}(N,\,D/N,\,\mathcal{A}_{N}/N).$$ Therefore, we have,
$$[\mathcal{K}(N,\,D/N,\,\mathcal{A}_{N}/N):\mathbb{F}_{q}] =
[\mathbb{F}_{q}(\zeta):\mathbb{F}_{q}]/[E_{G}(\mathcal{A}_{N}/D):\mathcal{A}_{N}]= o(\mathcal{A}_{N},\,D).$$  This 
proves  (\ref{E5}) and we thus have
\begin{eqnarray*}\mathbb{F}_{q}[G] &\cong& \bigoplus\limits_{d \mid [G:G']\atop
{ l \mid [\mathbb{F}_{q}(\xi):
\mathbb{F}_{q}]}} \bigoplus_{(N,\,D/N,\,\mathcal{A}_{N}/N)\in
\mathcal{S}_{d,\,l}}{M_{d}(\mathbb{F}_{q^{l}})}^{(|R(\mathcal{A}_{N}/D)|)}\\ \\
&\cong& \bigoplus\limits_{d \mid [G:G']\atop
{ l \mid [\mathbb{F}_{q}(\xi):
\mathbb{F}_{q}]}} {M_{d}(\mathbb{F}_{q^{l}})}^{(\alpha_{d,\,l})},\end{eqnarray*}
where $\alpha_{d,\,l} = \displaystyle\sum_{(N,\,D/N,\,\mathcal{A}_{N}/N)\in
\mathcal{S}_{d,\,l}}| R(\mathcal{A}_{N}/D)|.$  This proves $(i).$\para
$(ii)$~It follows from $(i)$ and the standard results on automorphisms of finite dimensional algebras.~~$\Box$
\section{  Metacyclic groups}
 We now  illustrate Theorem \ref{T1} with its application to metacyclic groups; thus obtaining an
alternative way of finding a complete set of primitive central idempotents of $\mathbb{F}_{q}[G]$ with $G$ given by 
presentation
\begin{equation} \label{e0} G = \langle a, b \,|\, a^{n} = 1,\, b^{t} = a^{k}, \, b^{-1}ab=
a^{r}\rangle,\end{equation} where  $n,\,t,\,k,\,r$ are  natural numbers with $r^{t}\equiv 1~({\rm mod}\,n),\,
k(r-1)\equiv 0~({\rm mod}\, n).$

 \par\vspace{.5cm} For a divisor $v$ of $n$, let\begin{itemize}
\item $o_{v} = {\rm {ord}}_{v}(r). $
 \item $G_{o_{v}} = \langle a,\,b^{o_{v}}\rangle$.
\item  $\mathcal{B}_{o_{v}} = \{(w,\,i,\,c)\in \mathbb{Z}^{3}\,|\,w\,>\,0,\, 
    w\,|\,n,\,w\,|\,r^{o_{v}}-1,\,o_{v}c\,>\,0\,,
o_{v}c\,|\,t,\,w\,|\,
k+i\frac{t}{o_{v}c}\}.$
\end{itemize}
Let
$$\mathfrak{N} = \{(v,\,i,\,c)\in \mathbb{Z}^{3}\,\mid \,
v>0 ,\,v|n,\,
c>0,\, c|t,\,0\leq i \leq v-1,\,v|k+i\frac{t}{c},\,o_{v}\,|\,c ~ \mbox{and}~  v\,|\,i(r-1)\}.$$
For $(v,\,i,\,c)\in \mathfrak{N},$ define \begin{itemize}
%\item
\item $H_{v,\,i,\,c} = \langle a^{v},\,a^{i}b^{c}\rangle.$
\item $X_{v,\,i,\,c} = \{(v,\,\alpha,\beta)\,|\,\beta o_{v}\mid c,\,\alpha \frac{c}{\beta o_{v}}\equiv i~ ({\rm 
    mod}\,v),\,
\beta = \frac{c \gcd(\alpha(r-1),\,v)}{v o_{v}},$\linebreak
\,gcd$(v,\,\alpha,\,\beta) =1$ ~\mbox{and}~$(v,\,\alpha,\,\beta)\in
\mathcal{B}_{o_{v}}\}.$
\end{itemize}\par\vspace{.25cm}
Define a relation, denoted $\sim,$ on $X_{v,\,i,\,c}$ as follows:\vspace{2mm}\\
For $(v,\,\alpha_{1},\,\beta_{1}),\,(v,\,\alpha_{2},\,\beta_{2})\in X_{v,\,i,\,c} ,$ we say that 
$(v,\,\alpha_{1},\,\beta_{1})\sim
(v,\,\alpha_{2},\,\beta_{2})\Leftrightarrow \beta_{1}=\beta_{2}$ and $\alpha_{1}\equiv \alpha_{2} r^{j}~({\rm 
mod}\,v)$
for some $j.$ It is easy
to see that $\sim$ is an equivalence relation on $X_{v,\,i,\,c}.$ Let 
$\mathfrak{X}_{v,\,i,\,c}$\index{$\mathfrak{X}_{v,\,i,\,c}$} denote
 the set of distinct equivalence classes of $X_{v,\,i,\,c}$ under the equivalence relation
$\sim.$\vspace{.15cm}

\begin{theorem}\label{T11} Let $\mathbb{F}_{q}$ be a finite field with $q$ elements and $G$ the group given by the 
presentation
(\ref{e0}). If $\gcd
(q,\,nt) = 1,$ then
$$\bigcup_{(v,\,i,\,c)\in \mathfrak{N}}\{e_{C}(G,\,G_{o_{v}},\,H_{v,\,\alpha,\,\beta o_{v}}),~\mid~
(v,\,\alpha,\,\beta)\in
\mathfrak{X}_{v,\,i,\,c},~C\in \mathcal{R}(G_{o_{v}}/H_{v,\,\alpha,\,\beta o_{v}}) \}$$ is a complete set of primitive 
central idempotents
of the group algebra  $\mathbb{F}_{q}[G].$
\end{theorem}\index{primitive central idempotents!metacyclic group algebras} \bigskip
We prove the above result in a number of steps.\bigskip

\begin{lemma}\label{L1}
$H_{v,\,i,\,c},~ (v,\,i,\,c)\in \mathfrak{N},$ are all the distinct normal subgroups of $G.$
 \end{lemma}\par\vspace{.25cm}
\noindent {\bf{Proof}.} Let $N\unlhd G.$ Suppose $N\cap \langle a \rangle = \langle a^{v} \rangle, v\,|\,n,\,v\,>\,0.$ 
Now, if
 $N/N\cap \langle a \rangle,$ as a subgroup of $G/\langle a \rangle,$ is generated by $\langle b^{c}\langle a
\rangle\rangle,\,c\,>\,0,\,c\,|\,t,$ then clearly,
 \begin{equation} \label{f0'}N = \langle a^{v},\,a^{i}b^{c} \rangle~\mbox{for some}~i,\,0\leq i \leq 
 v-1.\end{equation}
\par\vspace{.25cm}
Now $N$ being a normal subgroup of $G,$ we must have $ b^{-1}a^{i}b^{c}b,~ a^{-1}a^{i}b^{c}a$ and 
${(a^{i}b^{c})}^{t/c}$
all belong to $N.$ This gives
\begin{equation}\label{f0} v\,|\,i(r-1),~o_{v}\,|\,c,~ v\,|\,k+i\frac{t}{c}.\end{equation}
Consequently, equations (\ref{f0'}) and (\ref{f0}) yield that $(v,\,i,\,c)\in \mathfrak{N}$ and $N= H_{v,\,i,\,c}.$
\par\vspace{.25cm}
Conversely, it is easy to see that for any  $(v,\,i,\,c)\in \mathfrak{N},\, H_{v,\,i,\,c}$ is normal subgroup of $G.$ 
Furthermore,
\begin{equation}\label{f1}
 |H_{v,\,i,\,c}|=\frac{nt}{vc}.
\end{equation}

\par\vspace{.25cm} In order to complete the proof of the Lemma, we need to show that\linebreak $H_{v,\,i,\,c},$
$ (v,\,i,\,c)\in \mathfrak{N},$
 are distinct. Let
$(v_{1},\,i_{1},\,c_{1}),~(v_{2},\,i_{2},\,c_{2})\in \mathfrak{N}$ be such that
$H_{v_{1},\,i_{1},\,c_{1}} =
H_{v_{2},\,i_{2},\,c_{2}}.$ Then $ \langle a^{v_{1}} \rangle = H_{v_{1},\,i_{1},\,c_{1}}\cap \langle a \rangle =
H_{v_{2},\,i_{2},\,c_{2}}\cap \langle a \rangle =
\langle a^{v_{2}}\rangle  $\linebreak implies  that $ v_{1} = v_{2} = v,$ say. Also, in view of equation (\ref{f1}),
$|H_{v,\,i_{1},\,c_{1}}/\langle a^{v} \rangle| = |H_{v,\,i_{2},\,c_{2}}/\langle a^{v} \rangle|$ implies that
$c_{1} = c_{2}= c,$ say. Further,
  $a^{i_{2}}b^{c}\in H_{v,\,i_{2},\,c},$
 $a^{i_{1}}b^{c}\in H_{v,\,i_{1},\,c}$ and $H_{v,\,i_{1},\,c} = H_{v,\,i_{2},\,c}$ gives that $a^{i_{1}-i_{2}}\in 
 H_{v,\,i_{1},\,c}
\cap\langle a \rangle = \langle a^{v} \rangle.$ Hence $i_{2}\equiv
i_{1}~({\rm mod}\,v_{1}),$ i.e., $i_{1}=i_{2}.$ This proves the Lemma.  $\Box$\par\vspace{.5cm}

 \begin{lemma}\label{L5}~ Let $(v,\,i,\,c)\in \mathfrak{N}$ and $N = H_{v,\,i,\,c}.$ Then \vspace{2mm}\\
$(i)$~ $G_{o_{v}}/N$ is a maximal abelian subgroup
of $G/N$ containing ${(G/N)}^{'}.$\vspace{.1cm}\\
$(ii)$~$H/N$ is a subgroup of $G_{o_{v}}/N$ with
cyclic quotient and $H/N$ core-free in $G/N \Leftrightarrow H = H_{v,\,\alpha,\,\beta o_{v}},~(v,\,\alpha,\,\beta)
\in X_{v,\,i,\,c}.$
\end{lemma}\par\vspace{.25cm}
\noindent{\bf{Proof.}}~$(i)$~ By (\cite{cur}, p.336), $G_{o_{v}}'= \langle a^{r^{o_{v}}-1}\rangle.$ Since 
$v\,|\,r^{o_{v}}-1,$
 we have \linebreak
$G_{o_{v}}'\leq \langle a^{v} \rangle \leq N$ and therefore $G_{o_{v}}/N$ is abelian. Furthermore, $G_{o_{v}}/N$ 
contains
 ${(G/N)}^{'}$
as $G' = \langle a^{r-1} \rangle \leq \langle a,\,b^{o_{v}} \rangle = G_{o_{v}}.$ Thus $G_{o_{v}}/N$ is an abelian 
subgroup
of $G/N$ containing ${(G/N)}^{'}.$\par\vspace{.25cm}  If $o_{v}=1,$ then clearly, $G_{o_{v}}/N =G/N$ is
a maximal abelian subgroup
of $G/N$ containing ${(G/N)}^{'}.$ Let $o_{v}> 1.$
Suppose that $K/N$ is an abelian subgroup of $G/N$ with  $G_{o_{v}}/N\leq K/N \leq G/N.$ Since $o_{v}> 1,~ G/N$ is not 
abelian. Thus
$K/N\lneq G/N.$ Now  $ K\cap \langle a \rangle = \langle a \rangle$ implies that  $K = \langle a ,\,b^{j}\rangle $ for 
some $j\,|\,o_{v}.$
 However, $K' \leq N$ implies that $ \langle a^{r^{j}-1}\rangle\leq N,$ which gives that
 $v\,|\,r^{j}-1$, i.e.,  $o_{v}\,|\,j.$
 Thus $j=o_{v}$ and $K/N = G_{o_{v}}/N.$ This proves $(i).$
\par\vspace{.25cm}
$(ii)$~Let $H/N$ be a subgroup of $G_{o_{v}}/N$ with
cyclic quotient. By (\cite{oli1}, Lemma 2.2), we have $$H = H_{u,\,\alpha,\,\beta o_{v}},\,
(u,\,\alpha,\,\beta)\in \mathcal{B}_{o_{v}}~\mbox{and}~\gcd(u,\,\alpha,\,\beta) =1.$$
 Since $N\,\leq\,H$, we must have  $a^{v}\in H$ and $a^{i}b^{c}\in H,$ which holds, if, and only if,
\begin{equation}\label{m1}
u\,|\,v,~\beta o_{v}\,|\,c~\mbox{and}~ \alpha \frac{c}{\beta o_{v}}\equiv i~({\rm {mod}}\,u).
 \end{equation}\par\vspace{.25cm}
We claim that $\operatorname{core}(H),$ the largest normal subgroup of $G$ contained in $H,$ is given by
$$\operatorname{core}(H)= \langle a^{u},\,a^{\alpha\frac{\delta}{\beta o_{v}}}b^{\delta} \rangle,\,\delta =
\frac{\beta u o_{v}}{\gcd(\alpha(r-1),\,u)}.$$
Let $K = \langle a^{u},\,a^{\alpha \frac
 {\delta}{\beta o_{v}}}b^{\delta} \rangle$ with $\delta$ as above. Since $(u,\,\alpha\frac{\delta}{\beta 
 o_{v}},\,\delta)
\in \mathfrak{N},$ by  Lemma \ref{L1}, it follows that $K$ is a normal subgroup of $G.$ Since
$ab^{o_{v}}a^{-1}b^{-o_{v}} \in \langle a^{v}\rangle,$ we have
 $a^{\alpha \frac{\delta}{\beta o_{v}}} b^{\delta}
{ (a^{\alpha}b^{\beta o_{v}})}^ {-\frac{\delta}{\beta o_{v}}}\in \langle a^{v}\rangle.$ Thus $K$ is a subgroup of 
$H_{u,\,\alpha,\,\beta o_{v}} = H.$ \par\vspace{.25cm}
In order to show that $\operatorname
{core}(H) = K,$ we need to  show that $K$ is the largest normal subgroup of $G$ contained  in
$H= H_{u,\,\alpha,\,\beta o_{v}}.$ Let $L$ be a normal subgroup of $G$ contained in $H_{u,\,\alpha,\,\beta o_{v}}.$ By 
Lemma \ref{L1},
$L = H_{w,\,\gamma,\,f}$ for some $(w,\,\gamma,\,f) \in \mathfrak{N}.$ Since $\langle a^{w} \rangle  =
L \cap \langle a \rangle \leq H_{u,\,\alpha,\,\beta o_{v}}\cap \langle a \rangle =  \langle a^{u} \rangle ,$
 it follows that
$u\,|\, w.$ Next observe that an arbitray element of $ H_{u,\,\alpha,\,\beta o_{v}}$ is of the type $a^{j}b^{s}$ with
$\beta o_{v}\,|\,s$ and \linebreak$j \equiv \alpha\frac{s}{\beta o_{v}}~({\rm mod}\,u).$ Therefore, $L = 
H_{w,\,\gamma,\,f}$ is a subgroup of
$H_{u,\,\alpha,\,\beta o_{v}}$ if, and only if, $\beta o_{v}\,|\,f$ and $\gamma  \equiv \alpha\frac{f}{\beta o_{v}}~
({\rm mod}\,u).$ Since $\gamma(r-1)\equiv 0~({\rm mod}\,w),$
we have\linebreak $\alpha\frac{f}{\beta o_{v}}(r-1)\equiv 0~({\rm{mod}}\,u).$ This gives that $\delta\,|\,f$  and 
consequently
$L = H_{w,\,\gamma,\,f}$ is contained in $K = \langle a^{u},\,a^{\alpha \frac {\delta}{\beta o_{v}}} b^{\delta} 
\rangle.$ This proves that $K$ is the largest normal subgroup of $G$ contained in
$H_{u,\,\alpha,\,\beta o_{v}},$ which proves the claim.
\par\vspace{.25cm}
It is now immediate from the claim that  $H/N$ is core-free in $G/N$ if, and only if, $u=v$ and
$\delta = c.$ This proves $(ii).$  $\Box$
\par\vspace{.25cm}

\begin{lemma}\label{L7}~ Let $(v,\,i,\,c)\in \mathfrak{N}$ and $(v,\,\alpha_{1},\,\beta_{1}),~
(v,\,\alpha_{2},\,\beta_{2})\in X_{v,\,i,\,c}.$   Then $H_{v,\,\alpha_{1},\,\beta_{1}o_{v}}$ and
$H_{v,\,\alpha_{2},\,\beta_{2}o_{v}}$ are conjugate in $G$ if, and only if, $\beta_{1} = \beta_{2}$ and
$\alpha_{1}\equiv \alpha_{2}r^{j}~({\rm mod}\,v),$ for some $j.$
 \end{lemma}\vspace{.25cm}
{\bf{Proof.}} Suppose
\begin{equation}\label{m0'}
 H_{v,\,\alpha_{1},\,\beta_{1}o_{v}} = g^{-1}H_{v,\,\alpha_{2},\,\beta_{2}o_{v}} g,~~~g = a^{i}b^{j}\in 
 G.\end{equation}
 Then, in particular, in view of  equation (\ref{f1}), we have
$$|H_{v,\,\alpha_{1},\,\beta_{1}o_{v}}| = \frac{nt}{v\beta_{1}o_{v}} = \frac{nt}{v\beta_{2}o_{v}}=
|H_{v,\,\alpha_{2},\,\beta_{2}o_{v}}|,$$ i.e.,
$$\beta_{1} = \beta_{2}.$$\par\vspace{.25cm}
Further equation (\ref{m0'}) holds, if, and only if,  
$${(a^{i}b^{j})}^{-1}a^{\alpha_{2}}b^{\beta_{1}o_{v}}a^{i}b^{j}\in
 H_{v,\,\alpha_{1},\,\beta_{1}o_{v}}.$$
 Since $ab^{o_{v}}a^{-1}b^{-o_{v}}\in \langle a^{v}\rangle,$ we have
${(a^{i}b^{j})}^{-1}a^{\alpha_{2}}b^{\beta_{1}o_{v}}a^{i}b^{j} {(a^{\alpha_{2}r^{j}}b^{\beta_{1}o_{v}})}^{-1}\in
\langle a^{v}\rangle\subseteq H_{v,\,\alpha_{1},\,\beta_{1}o_{v}},$ which yields that
\begin{equation*}\alpha_{1} \equiv
\alpha_{2}r^{j}~({\rm {mod}}\, v)\end{equation*}  and  proves the Lemma.  $\Box$
\vspace{.5cm}\\
{\bf{Proof of Theorem \ref{T11}.}}~ By Lemma \ref{L1}, $H_{v,\,i,\,c},\,
(v,\,i,\,c)\in \mathfrak{N},$ are all the distinct normal subgroups of $G.$ For $(v,\,i,\,c)\in \mathfrak{N},$ and
$N=H_{v,\,i,\,c},$ Lemma \ref{L5}  implies that
$$ S_{G/N} = \{(H_{v,\,\alpha,\,\beta o_{v}}/N,\,G_{o_{v}}/N)\,|\,(v,\,\alpha,\,\beta)\in
\mathfrak{X}_{v,\,i,\,c} \}.$$
Therefore, we have
$$\mathcal{S} = \bigcup_{(v,\,i,\,c)\in \mathfrak{N}} \{(H_{v,\,i,\,c},\,H_{v,\,\alpha,\,\beta 
o_{v}}/N,\,G_{o_{v}}/N)\,|\,(v,\,\alpha,\,\beta)\in
\mathfrak{X}_{v,\,i,\,c} \}$$ and consequently, Theorem \ref{T1} yields the required result.  $\Box$
{\section{Groups with central quotient Klein four-group}}
The groups $G$ of the type $G/Z(G) \cong \mathbb{Z}_{2}\times \mathbb{Z}_{2},$ where $Z(G)$ denotes the centre of
the group $G$, arose in the work of Goodaire \cite{goo83} while studying Moufang loops. It is known (\cite{goo1}, 
Chapter 5) that any group with $G/Z(G)\cong\mathbb{Z}_{2}\times \mathbb{Z}_{2}$ is the direct product of an 
indecomposable
group (with this property) and an abelian group; Moreover the finite indecomposable groups with $G/Z(G)\cong 
\mathbb{Z}_{2}\times \mathbb{Z}_{2}$
break into five classes as follows:
\par\vspace{.5cm}
\begin{tabular}{|c|c|c|}\hline  Group&Generators&Relations\\ \hline $D_{1}$&$x,y,t$&$x^2,
y^2,t^{2^{m}},y^{-1}x^{-1}yxt^{2^{m-1}},t~\mbox{central}~,$ \\ &&$m
\geq 1$ \\ \hline $D_{2}$&$x,y,t$&$x^2t^{-1},
y^2t^{-1},t^{2^{m}},y^{-1}x^{-1}yxt^{2^{m-1}},t~\mbox{central}~,$ \\ &&$m
\geq 1$\\ \hline $D_{3}$&$x,y,t_{1},t_{2}$&$x^2,
y^2t_{2}^{-1},t_{1}^{2^{m_{1}}},t_{2}^{2^{m_{2}}},y^{-1}x^{-1}yxt_{1}^{2^{m_{1}-1}},t_{1},t_{2}~\mbox{central}~,$
 \\ &&$m_{1},m_{2}
\geq 1$\\ \hline $D_{4}$&$x,y,t_{1},t_{2}$&$x^2t_{1}^{-1}
,
y^{2}t_{2}^{-1},t_{1}^{2^{m_{1}}},t_{2}^{2^{m_{2}}},y^{-1}x^{-1}yxt_{1}^{2^{m_{1}-1}},t_{1},t_{2}~\mbox{central}~$,
\\ &&$m_{1},m_{2}
\geq 1$\\ \hline $D_{5}$&$x,y,t_{1},t_{2},t_{3}
$&$x^{2}t_{2}^{-1},
y^{2}t_{3}^{-1},t_{1}^{2^{m_{1}}},t_{2}^{2^{m_{2}}},t_{3}^{2^{m_{3}}},y^{-1}x^{-1}yxt_{1}^{2^{m_{1}-1}},
t_{1},t_{2},t_{3}~\mbox{central}~,$\\ &&$m_{1},m_{2},m_{3}
\geq 1$  \\ \hline\end{tabular}\par\vspace{.5cm}
It thus becomes important to investigate the group algebra $\mathbb{F}_{q}[G],$ $G$ of type $D_{i}$, $1\leq i \leq 5.$ 
Recently, Ferraz, Goodaire and Milies \cite{fer2} have given a lower bound on the number of simple components of
these semisimple finite group algebras. We improve  Theorems 3.1 \& 3.2 of
 \cite{fer2} by providing  the complete algebraic structure of $\mathbb{F}_{q}[G]$,   $G$ of type $D_{i}$, $i = 1, 2$.

\bigskip
{\subsection{Groups of type $D_{1}$}}Observe that for $m=1,$ the group $G$ of type $D_{1}$ is isomorphic to $D_{8},$ the dihedral group of
order $8,$ and the structure of
group algebra $\mathbb{F}_{q}[D_{8}]$ can be read from Example $4.3$ of \cite{sha2}.\par\vspace{.25cm}
Let $m\geq 2,$ Define\\\indent
 $N_{0}:=\langle 1 \rangle,~N_{1}:=\langle t,\,x\rangle,~N_{2}:=\langle t,\,y\rangle,~
N_{3}:=\langle t,\,xy\rangle,
N_{4}^{(\alpha)}:=\langle t^{2^{\alpha}},\,x,\,y\rangle,$\\\indent $ N_{5}^{(\beta)}:=\langle 
t^{2^{m-1}},\,x,\,yt^{2^{\beta}}
\rangle,
N_{6}^{(\beta)}:=\langle t^{2^{m-1}},\,xt^{2^{\beta}},\,y\rangle,
N_{7}^{(\beta)}:=\langle t^{2^{\beta}}x,
t^{2^{\beta}}y\rangle,$\\\indent $0\leq \alpha \leq m-1,\,0\leq \beta \leq m-2.$\par
Let $\lambda$ be the highest power of $2$ dividing $q-1$
(resp. $q+1$) according as $q\equiv 1~({\rm mod}\,4)$ (resp. $q\equiv -1~({\rm mod}\,4)$).\par\vspace{.25cm}
 Ferraz, Goodaire and Milies (\cite{fer2}, Theorem 3.1) proved
 that the Wedderburn decomposition of $\mathbb{F}_{q}[G],\,G$ of type $D_{1},$ contains at least $8m-10$
simple components. If $q\equiv 3 ~({\rm mod}\,8),$ then this number is acheived with $8m-12$ fields and
$2$ quaternion algebras, each necessarily a ring of $2\times 2$ matrices. The following result improves the result of  
Ferraz et.al. by providing a concrete description of $\mathbb{F}_{q}[G],~G$ of type $D_{1}$:
 \pagebreak
\begin{theorem} A complete set of primitive  central idempotents, Wedderburn decomposition and the automorphism group 
of
$\mathbb{F}_{q}[G],\,G$ of type $D_{1},\, m\geq 2$, are given by:
\par\vspace{.25cm}
   \begin{center} {\bf{Primitive central
idempotents}} \index{primitive central idempotents of $\mathbb{F}_{q}[G]$!$G$ of type $D_{1}$}\end{center}
\vspace{.25cm}\begin{quote}
$e_{C}(G,\,N_{1},\,\langle x \rangle),\,C\in R(N_{1}/\langle x \rangle);\\
e_{C}(G,\,G,\,N_{i}),\,C\in R(G/N_{i}),\,1\leq i\leq 3;\\
e_{C}(G,\,G,\,N_{4}^{(\alpha)}),\,C\in R(G/N_{4}^{(\alpha)}),\,0\leq \alpha \leq m-1;\\
e_{C}(G,\,G,\,N_{j}^{(\beta)}),\,C\in R(G/N_{j}^{(\beta)}),\,0\leq \beta \leq m-2,\,5\leq j\leq 7.$\end{quote}
\par\vspace{.25cm}
 \begin{center}{\bf{Wedderburn decomposition}}
\index{Wedderburn decomposition of $\mathbb{F}_{q}[G]$!$G$ of type $D_{1}$}\end{center}

\underline{$q\equiv 1 ~({\rm mod}\,4)$}

$$\mathbb{F}_{q}[G]\cong \begin{cases}{\mathbb{F}_{q}}^{(2^{m+1})}\bigoplus {M_{2}(\mathbb{F}_{q})}^{(2^{m-1})},&
m\leq \lambda,\\{\mathbb{F}_{q}}^{(2^{m+1})}\bigoplus {M_{2}(\mathbb{F}_{q^2})}^{(2^{m-2})},&
m= \lambda+1,\\ {\mathbb{F}_{q}}^{(2^{\lambda+2})}\bigoplus_
{\alpha =\lambda+1}^{m-1} {\mathbb{F}_{q^{2^{\alpha - \lambda}}}}^{(2^{\lambda +1})}
 \bigoplus {M_{2}(\mathbb{F}_{q^{2^{m-\lambda}}})}^{(2^{\lambda-1})},&
m\geq \lambda+2.\end{cases}$$\par\vspace{.25cm}
\underline{$q \equiv -1~({\rm mod}\,4)$}
$$\mathbb{F}_{q}[G]\cong \begin{cases}{\mathbb{F}_{q}}^{(8)} \bigoplus
{\mathbb{F}_{q^{2}}}^{(2^{m}-4)}\bigoplus {M_{2}(\mathbb{F}_{q^{2}})}^{(2^{m-2})},&
2\leq m\leq \lambda+1,\\ {\mathbb{F}_{q}}^{(8)} \bigoplus
{\mathbb{F}_{q^{2}}}^{(2^{m}-4)}\bigoplus {M_{2}(\mathbb{F}_{q^{4}})}^{(2^{m-3})},&
 m=\lambda+2,\\  {\mathbb{F}_{q}}^{(8)} \bigoplus{\mathbb{F}_{q^{2}}}^{(2^{\lambda+2}-4)}\bigoplus_
{\alpha =\lambda+2}^{m-1} {\mathbb{F}_{q^{2^{\alpha - \lambda}}}}^{(2^{\lambda+1})}
 \bigoplus {M_{2}(\mathbb{F}_{q^{2^{m-\lambda}}})}^{(2^{\lambda-1})},&
m\geq \lambda+3.\end{cases}$$
\par\vspace{.25cm}
 \begin{center}{\bf{Automorphism group}}\end{center}\underline{$q\equiv 1 ~({\rm mod}\,4)$}
$$\mathrm{Aut}(\mathbb{F}_{q}[G])\cong \begin{cases}S_{2^{m+1}}\bigoplus
({\mathrm{SL}_{2}(\mathbb{F}_{q})}^{(2^{m-1})}\rtimes S_{2^{m-1}}),&
 m\leq \lambda ,\\ S_{2^{m+1}}\bigoplus
\left(({\mathrm{SL}_{2}(\mathbb{F}_{q^2})\rtimes \mathbb{Z}_{2})}^{(2^{m-2})}\rtimes S_{2^{m-2}}\right),&
 m= \lambda+1,\\ S_{2^{\lambda+2}}\bigoplus_{\alpha=\lambda+1}^{m-1}
({\mathbb{Z}_{2^{\alpha-\lambda}}}^{(2^{\lambda+1})}\rtimes S _{2^{\lambda+1}})\bigoplus
\mathcal{H}_{\lambda},& m\geq \lambda+2,\end{cases}$$ \par\vspace{.25cm}\pagebreak
\underline{$q\equiv -1 ~({\rm mod}\,4)$}
$$\mathrm{Aut}(\mathbb{F}_{q}[G])\cong \begin{cases}S_{8}\bigoplus ({\mathbb{Z}_{2}}^{(2^{m}-4)}\rtimes S_{2^{m}-4})
\bigoplus(({\mathrm{SL}_{2}(\mathbb{F}_{q^{2}})\rtimes \mathbb{Z}_{2})}^{(2^{m-2})}
\rtimes S_{2^{m-2}}),&
 m\leq \lambda+1,\\
S_{8}\bigoplus ({\mathbb{Z}_{2}}^{(2^{m}-4)}\rtimes S_{2^{m}-4})
\bigoplus(({\mathrm{SL}_{2}(\mathbb{F}_{q^{4}})\rtimes \mathbb{Z}_{4})}^{(2^{m-3})}
\rtimes S_{2^{m-3}}),&
 m= \lambda+2,\\ S_{8}\bigoplus({\mathbb{Z}_{2}}^{(2^{\lambda+2}-4)}\rtimes S_{2^{\lambda+2}-4})
 \bigoplus_{\alpha=\lambda+2}^{m-1}({\mathbb{Z}_{2^{\alpha-\lambda}}}^{(2^{\lambda+1})}\rtimes S 
 _{2^{\lambda+1}})\bigoplus
\mathcal{H}_{\lambda},& m\geq \lambda+3,\end{cases}$$ where $\mathcal{H}_{\lambda} = {(\mathrm{SL}_{2}
(\mathbb{F}_{q^{2^{m-\lambda}}})\rtimes \mathbb{Z}_{2^{m-\lambda}})}^{(2^{\lambda-1})}\rtimes S_{2^{\lambda-1}}.$
\end{theorem}\bigskip
In order to prove the above Theorem, we first need to compute all the normal subgroups of $G,~G$ of type 
$D_{1}.$\bigskip
\begin{lemma}\label{l.5}All the distinct non-identity normal subgroups of $G$ are given by:

 \begin{quote}$(i)$~$\langle t^{2^{\alpha}},\,x \rangle,\,\langle t^{2^{\alpha}},\,y \rangle,\,
\langle t^{2^{\alpha}},\,xy \rangle,\,\langle t^{2^{\alpha}},\,x,\,y \rangle;$\vspace{2mm}\\
$(ii)$~$\langle t^{2^{\beta}}x \rangle,\,\langle t^{2^{\beta}}y \rangle,\,\langle t^{2^{m-1}},\,t^{2^{\beta}}xy 
\rangle,\,
\langle t^{2^{m-1}},\,x,\,t^{2^{\beta}}y \rangle,\,\langle t^{2^{m-1}},\,t^{2^{\beta}}x,\,y \rangle,
\,\langle t^{2^{\beta}}x,\,t^{2^{\beta}}y \rangle;$\vspace{2mm}\\
$(iii)$~$\langle t^{2^{\gamma}}\rangle,$\end{quote} where $0\leq \alpha\leq m-1,~
 0\leq \beta\leq m-2$ and $0\leq \gamma\leq m-1.$
\end{lemma}\par\vspace{.25cm}
\noindent{\bf{Proof.}} Observe that all the subgroups listed in the statement are distinct and normal in 
$G.$\par\vspace{.25cm}

Let $N$ be a normal subgroup of $G$ not contained in $\langle t \rangle.$ If $N\neq \langle 1 \rangle,$ then it is easy 
to see that
 $\langle t^{2^{m-1}}
 \rangle \leq N.$ Therefore $N\cap \langle t \rangle =
\langle t^{2^{v}} \rangle ,\,0\leq v\leq m-1.$ Since
$N/N\cap \langle t \rangle$ is isomorphic to subgroup of $G/\langle t \rangle,$ which is generated by
$x\langle t \rangle,\, y\langle t \rangle,$ it follows that
$N/N\cap \langle t \rangle$ is isomorphic to one of the following:
$ \langle x\langle t \rangle\rangle,\,\langle y\langle t \rangle\rangle,\,\langle xy\langle t \rangle\rangle$
or $\langle x\langle t \rangle,\,y\langle t \rangle\rangle.$\par\vspace{.25cm}
\noindent {\bf{Case I :}}~ $N/ \langle t^{2^{v}} \rangle
\cong \langle x\langle t \rangle\rangle$\vspace{.25cm}\\ In this case,  $N =
\langle t^{2^{v}},\,t^{2^{i}}x\rangle,$ for some $i,\, 0\leq i \leq v\leq m-1.$\par\vspace{.15cm}
 If $i=v,$ then $N =
\langle t^{2^{v}},\,x\rangle.$ Since $N \unlhd G,\,xt^{2^{m-1}}=y^{-1}xy\in N,$ implies that $t^{2^{m-1}}
\in N\cap \langle t \rangle = \langle t^{2^{v}} \rangle,$ which is possible only if $v\leq m-1.$ \par\vspace{.15cm}
 If $i<v,$ then
 $N =
\langle t^{2^{v}},\,t^{2^{i}}x\rangle =
\langle t^{2^{i}}x\rangle$ as $t^{2^{v}}\in \langle t^{2^{i}}x\rangle.$ Further $xt^{2^{i}+2^{m-1}}=y^{-1}t^{2^{i}}xy
\in N$ implies that $t^{2^{m-1}}\in \langle t^{2^{v}}\rangle.$ Hence $v\leq m-1$ and  $i\leq m-2.$\par\pagebreak Thus 
in this case, either
\begin{equation}
\label{e5.0} N =
\langle t^{2^{i}},\,x\rangle,\,0\leq i \leq m-1 \end{equation} or
\begin{equation}\label{e5.1}
 N = \langle t^{2^{i}}x\rangle,\,0\leq i \leq m-2.\end{equation} {\bf{Case II:}}~   $N/ \langle t^{2^{v}} \rangle
\cong \langle y\langle t \rangle\rangle.$\vspace{.25cm}\\
Computation analogous to those in Case I yield that
 \begin{equation}\label{e5.2} N =
 \langle t^{2^{i}}y\rangle,0\leq i\leq m-2\end{equation} or \begin{equation}\label{e5.3}
N = \langle  t^{2^{i}},\,y\rangle,\,0\leq i \leq m-1.
\end{equation}\par\vspace{.25cm}\noindent {\bf{Case III:}}~  $N/ \langle t^{2^{v}} \rangle
\cong \langle xy\langle t \rangle\rangle.$\vspace{.25cm}\\
In this case $N  = \langle t^{2^{v}},\,t^{2^{i}}xy \rangle$ for  $0\leq i \leq v \leq m-1.$
\par\vspace{.15cm}
If $i=v,$ then $N= \langle t^{2^{v}},\,xy \rangle.$ Since $N$ is a normal subgroup of $G,\, xyt^{2^{m-1}} = y^{-1}xyy
\in N,$ implies that $t^{2^{m-1}}\in N\cap \langle t \rangle = \langle t^{2^{v}} \rangle,$ which is possible only
if $v\leq m-1.$\par\vspace{.15cm}
 If $i<v,$ then $N  =\langle  t^{2^{v}},\,t^{2^{i}}xy\rangle,\,0\leq i\leq m-2.$
% \par\vspace{.25cm}
% Now for $0 \leq v \leq m-1,\, 0\leq i \leq m-2,$
 Since $\langle  t^{2^{m-1}},\,t^{2^{i}}xy\rangle \leq \langle  t^{2^{v}},\,t^{2^{i}}xy\rangle$ and
$$t^{2^{v}}=\begin{cases} {(t^{2^{i}}xy)}^{2^{v-i}}t^{2^{m-1}},&~\mbox{if}~v-i=1,\\
{(t^{2^{i}}xy)}^{2^{v-i}},&~\mbox{if}~v-i\geq 2,\end{cases}$$ it follows that
$\langle  t^{2^{v}},\,t^{2^{i}}xy\rangle =\langle  t^{2^{m-1}},\,t^{2^{i}}xy\rangle.$\vspace{.25cm}\\
Thus in this case, either  \begin{equation}\label{e5.4} N =
\langle  t^{2^{i}},\,xy\rangle,\,0\leq i \leq m-1 \end{equation} or
  \begin{equation}\label{e5.5} N =
\langle  t^{2^{m-1}},\,t^{2^{i}}xy\rangle,\,0\leq i \leq m-2. \end{equation}\par\vspace{.25cm}
\noindent {\bf{Case IV:}}~ $N/ \langle t^{2^{v}} \rangle
\cong \langle x\langle t \rangle,\,\langle y\langle t \rangle\rangle.$\vspace{.15cm}\\
 In this case, $N$ is one of the following forms:
\begin{quote}$(a)$~$ \langle t^{2^{v}},\,x,\,  y \rangle;$\\\pagebreak
$(b)$~$ \langle t^{2^{v}},\,t^{2^{i}}x,\,y \rangle$ for some $i,\,0\leq i \leq v-1;$\\
$(c)$~$ \langle t^{2^{v}},\,x,\,t^{2^{i}}y \rangle$ for some $i,\,0\leq i \leq v-1;$\\
$(d)$~ $\langle t^{2^{v}},\,t^{2^{i}}x,\,t^{2^{i}}y \rangle$ for some $i,\,0\leq i \leq v-1;$\\
$(e)$~$\langle t^{2^{v}},\,t^{2^{i}}x,\,t^{2^{j}}y \rangle$ for some $1\leq i,\,j \leq v-1,\,i\neq j.$
\end{quote}\par\vspace{.25cm}
 Observe that for $0 \leq i\leq v-1,$
$$\langle t^{2^{v}},\,t^{2^{i}}x,\,y \rangle = \langle t^{2^{m-1}},\,t^{2^{i}}x,\,y \rangle,$$
$$\langle t^{2^{v}},\,x,\,t^{2^{i}}y \rangle = \langle t^{2^{m-1}},\,x,\,t^{2^{i}}y \rangle,$$
 and
$$\langle t^{2^{v}},\,t^{2^{i}}x,\,t^{2^{i}}y \rangle = \langle t^{2^{m-1}},\,t^{2^{i}}x,\,t^{2^{i}}y \rangle.$$
Also for $1\leq i,\,j \leq v-1,\,i \neq j,$
$$ \langle t^{2^{v}},\,t^{2^{i}}x,\,t^{2^{j}}y \rangle = \begin{cases}\langle t^{2^{i}}x,\,y \rangle,& 
~\mbox{if}~i<j,\\
\langle x,\,t^{2^{j}}y \rangle,& ~\mbox{if}~j<i. \end{cases}$$ \par\vspace{.25cm}
Thus we have proved that any normal subgroup of $G$ not contained in $\langle t \rangle$ is one of the forms given in 
$(i)$ and $(ii)$
of the statement. This proves the Lemma.  $\Box$
\vspace{.5cm}\\{\bf Proof of Theorem 5.}
In order to apply Theorem \ref{T1} to a group $G$ of type $D_{1},$ we need to compute $S_{G/N}$ for all
normal subgroups $N$ of $G$  given by Lemma \ref{l.5}. \par\vspace{.25cm}
 Clearly if
$N=\langle 1 \rangle,\,S_{G/N}=\{(\langle x \rangle,\,\langle  t,\,x \rangle)\}.$\par\vspace{.25cm}
Suppose $N$ is a non-identity normal subgroup of $G,$ then $N$ is  one of the subgroups listed in Lemma \ref{l.5}. 
Since
 $G'= \langle t^{2^{m-1}} \rangle \leq N,$ we have  $\mathcal{A}_{N}/N = G/N$ and the corresponding
$$ S_{G/N} = \begin{cases}\{(\langle 1\rangle,\,G/N)\},& \mbox{if}~G/N~\mbox{ is
cyclic},\\ \o{},& ~\mbox{otherwise}.\end{cases}$$
Next we see that among all the normal subgroups $N$ of $G$ stated in Lemma \ref{l.5}, only the following
subgroups $N$  satisfy the condition that
$G/N$ is cyclic;
$$N_{i},\,N_{4}^{(\alpha)},\,N_{j}^{(\beta)},\,1\leq i \leq 3,\,0\leq \alpha \leq m-1,\,5\leq j \leq 7,\,
0\leq \beta \leq m-2.$$ Therefore $\mathcal{S} = \{ (N_{0},\,\langle x \rangle,\,N_{1})\} \cup \{(N_{i},\,\langle 1 
\rangle,\,G/N_{i})\,
|\, 1 \leq i \leq 3\} \cup \{ (N_{4}^{(\alpha)},\,\langle 1 \rangle,\,G/N_{4}^
{(\alpha)})|$ \linebreak $ 0\leq \alpha \leq m-1\} \cup \{(N_{j}^{(\beta)},\,
\langle 1 \rangle,\,G/N_{j}^{(\beta)})\, | \,0\leq \beta \leq m-2\} $
and thus $(i)$ follows.
\par\vspace{.45cm} In order to prove $(ii)$ and $(iii)$, we first note that for any integer
 $\gamma \geq 2,$\par\pagebreak
$$ {\rm
ord}_{2^{\gamma}}(q) = \begin{cases} 2^{\gamma-\lambda},&  \gamma \geq \lambda+1,~q\equiv 
1~\mbox{or}~-1~(\rm{mod}\,4),\\
 1,&   \gamma \leq \lambda~,q\equiv 1~(\rm{mod}\,4),\\ 2,&  \gamma \leq \lambda,~q\equiv 
 -1~(\rm{mod}\,4).\end{cases}$$
  Direct calculations yield
that for each $(N, D/N, A_{N}/N) \in \mathcal{S}$, the corresponding $o(\mathcal{A}_{N},\,D)$
 and $|R(\mathcal{A}_{N}/D)|$  are as given by the following tables:
 \par\vspace{.5cm}
\noindent
{\bf Case I\,: $q\equiv 1~({\rm mod}\,4)$.}\\
 \par\vspace{.2cm}
\begin{tabular}{|c|c|c|}\hline ($N, \, D/N,\,\mathcal{A}_{N}/N)$& $o(\mathcal{A}_{N},\,D)$&$|R(\mathcal{A}_{N}/D)|$ \\ 
\hline $(N_{i},\,\langle 1 \rangle,\,G/N_{i}),$&1
&1 \\ $1\leq i\leq 3$&&\\ \hline
$(N_{4}^{(0)},\,\langle 1 \rangle,\,G/N_{4}^
{(0)}),$& $1$&
$1$ \\ \hline
$(N_{4}^{(\alpha)},\,\langle 1 \rangle,\,G/N_{4}^
{(\alpha)}),$& $\begin{cases}2^{\alpha-\lambda},&\alpha\geq \lambda+1,\\1,&\alpha\leq \lambda\end{cases}$&
$\begin{cases}2^{\lambda-1},&\alpha\geq \lambda+1,\\2^{\alpha-1},&\alpha\leq \lambda\end{cases}$ \\ $1\leq \alpha\leq 
m-1$&&
\\ \hline $(N_{j}^{(\beta)},\,
\langle 1 \rangle,\,G/N_{j}^{(\beta)}),$& $\begin{cases}2^{\beta+1-\lambda},&\beta\geq \lambda,\\1,&\beta\leq
\lambda-1\end{cases}$
&$\begin{cases}2^{\lambda-1},&\beta\geq \lambda,\\2^{\beta},&\beta\leq \lambda-1\end{cases}$ \\ $5\leq j \leq 7,\,
0\leq \beta \leq m-2$&& \\ \hline
$(N_{0},\, \langle x \rangle,\,N_{1}) $&$\begin{cases}2^{m-\lambda},&m\geq \lambda+1,\\1,&m\leq \lambda\end{cases}$&
$\begin{cases}2^{\lambda-1},&m\geq \lambda+1,\\2^{m-1},&m\leq \lambda\end{cases}$ \\ \hline
\end{tabular}
\par\vspace{.5cm} \noindent
{\bf Case II\,: $q \equiv -1~(\rm{mod}\,4).$}
 \par\vspace{.5cm}
\begin{tabular}{|c|c|c|}\hline ($N, \, D/N,\,\mathcal{A}_{N}/N)$& $o(\mathcal{A}_{N},\,D)$&$|R(\mathcal{A}_{N}/D)|$ \\ 
\hline $(N_{i},\,\langle 1 \rangle,\,G/N_{i}),$&1
&1 \\ $1\leq i\leq 3$&&\\ \hline
$(N_{4}^{(\alpha)},\,\langle 1 \rangle,\,G/N_{4}^
{(\alpha)}),$& $1$&
$1$ \\  $0\leq \alpha\leq 1$&&\\ \hline
$(N_{4}^{(\alpha)},\,\langle 1 \rangle,\,G/N_{4}^
{(\alpha)}),$& $\begin{cases}2^{\alpha-\lambda},&\alpha\geq \lambda+2,\\2,&\alpha\leq \lambda +1\end{cases}$&
$\begin{cases}2^{\lambda-1},&\alpha\geq \lambda+2,\\2^{\alpha-2},&\alpha\leq \lambda +1\end{cases}$ \\ $2\leq \alpha
\leq m-1$&&
\\ \hline $(N_{j}^{(0)},\,
\langle 1 \rangle,\,G/N_{j}^{(0)}),$& $1$
&$1$ \\ $5\leq j \leq 7,\,$&&
\\ \hline $(N_{j}^{(\beta)},\,
\langle 1 \rangle,\,G/N_{j}^{(\beta)}),$& $\begin{cases}2^{\beta+1-\lambda},&\beta\geq \lambda+1,\\2,&\beta\leq 
\lambda\end{cases}$
&$\begin{cases}2^{\lambda-1},&\beta\geq \lambda +1,\\2^{\beta-1},&\beta\leq \lambda\end{cases}$ \\ $5\leq j \leq 7,\,
1\leq \beta \leq m-2$&& \\ \hline
$(N_{0},\, \langle x \rangle,\,N_{1}) $&$\begin{cases}2^{m-\lambda},&m\geq \lambda+1,\\2,&m\leq \lambda\end{cases}$&
$\begin{cases}2^{\lambda-1},&m\geq \lambda+2,\\2^{m-2},&m\leq \lambda+1\end{cases}$ \\ \hline
\end{tabular}\par\vspace{.5cm} Thus, Theorem \ref{T2} with the help of above two tables yield $(ii)$ and $(iii).$
{\subsection{Groups of type $D_{2}$}} Observe that for $m=1,$ the group $G$ of type $D_{2}$ is isomorphic to $Q_{8},$ the
quaternion group of order $8$ and the structure of group algebra $\mathbb{F}_{q}[Q_{8}]$ can be read from
Example $4.4$ of \cite{sha2}.\par\vspace{.25cm}
Let $m\geq 2.$ Define\\\indent
 $K_{0}:=\langle 1 \rangle,\,K_{1}:=\langle x \rangle;
{K_{2}}^{(\alpha)}:= \langle x^{2^{\alpha}},\,x^{2^{\alpha}-1}y\rangle,\,
{K_{3}}^{(\beta)}:= \langle x^{2^{\beta}},\,x^{2^{\beta-1}-1}y\rangle,$\\\indent
$0\leq \alpha \leq m,\,1\leq \beta \leq m.$\par Let $\lambda$ be the highest power of $2$ dividing $q-1$
(resp. $q+1$) according as $q\equiv 1~({\rm mod}\,4)$ (resp. $q\equiv -1~({\rm mod}\,4)$).\par\vspace{.25cm}
Ferraz, Goodaire and Milies proved (\cite{fer2}, Theorem 3.2) that the Wedderburn decomposition of 
$\mathbb{F}_{q}[G],\,G$
of type $D_{2},$ contains at least $4m$ simple components. If $q\equiv 3 ~({\rm mod}\,8),$ then this number is acheived 
with
$4m-2$ fields and $2$ quaternion algebras, each necessarily a ring of $2\times 2$ matrices.
The following Theorem improves this result of Ferraz et.al.\bigskip
\begin{theorem} A complete set of primitive  central idempotents, Wedderburn decomposition and the automorphism group 
of
$\mathbb{F}_{q}[G],\,G$ of type $D_{2},\,m\geq 2,$ are given by:
\par\vspace{.25cm}
 \begin{center}{\bf{Primitive central
idempotents}}\end{center}
              \begin{quote} $e_{C}(G,\,K_{1},\,K_{0}),\,C\in R(K_{1}/K_{0});\\
e_{C}(G,\,G,\,K_{1}),\,C\in R(G/K_{1});\\
e_{C}(G,\,G,\,K_{2}^{(\alpha)}),\,C\in R(G/{K_{2}}^{(\alpha)}),\,0\leq \alpha\leq m;\\
e_{C}(G,\,G,\,K_{3}^{(\beta)}),\,C\in R(G/K_{3}^{(\beta)}),\,1\leq \beta \leq m.$\end{quote}
                                                                                             \par\vspace{.25cm}
\begin{center}{\bf{Wedderburn decomposition}}\end{center}\par\vspace{.25cm}

\underline{$q \equiv 1~({\rm mod}\,4)$}
$$\mathbb{F}_{q}[G]\cong \begin{cases}{\mathbb{F}_{q}}^{(2^{m+1})}\bigoplus {M_{2}(\mathbb{F}_{q})}^{(2^{m-1})},&
m\leq \lambda,\\ {\mathbb{F}_{q}}^{(2^{\lambda+1})}\bigoplus_
{\alpha =\lambda+1}^{m} {\mathbb{F}_{q^{2^{\alpha - \lambda}}}}^{(2^{\lambda})}
 \bigoplus {M_{2}(\mathbb{F}_{q^{2^{m-\lambda}}})}^{(2^{\lambda-1})},&
m\geq \lambda +1.\end{cases}$$
\par\vspace{.25cm}\pagebreak
\underline{$q \equiv -1~({\rm mod}\,4)$}
$$\mathbb{F}_{q}[G]\cong \begin{cases}{\mathbb{F}_{q}}^{(4)} \bigoplus
{\mathbb{F}_{q^{2}}}^{(2^{m}-2)}\bigoplus {M_{2}(\mathbb{F}_{q^{2}})}^{(2^{m-2})},&
2\leq m\leq \lambda +1,\\ {\mathbb{F}_{q}}^{(4)} \bigoplus{\mathbb{F}_{q^{2}}}^{(2^{\lambda+1}-2)}\bigoplus_
{\alpha =\lambda+2}^{m} {\mathbb{F}_{q^{2^{\alpha - \lambda}}}}^{(2^{\lambda})}
 \bigoplus {M_{2}(\mathbb{F}_{q^{2^{m-\lambda}}})}^{(2^{\lambda-1})},&
m\geq \lambda+2.\end{cases}$$
\par\vspace{.25cm}

\begin{center}{\bf{Automorphism group}}\end{center}\par\vspace{.25cm}
\underline{$q\equiv 1~({\rm mod}\,4)$}
$$\mathrm{Aut}(\mathbb{F}_{q}[G])\cong \begin{cases}S_{2^{m+1}}\bigoplus({\mathrm{SL}_{2}(\mathbb{F}_{q})}^{(2^{m-1})}
\rtimes S_{2^{m-1}}),&
 m\leq \lambda,\\ 
 S_{2^{\lambda+1}}\bigoplus_{\alpha=\lambda+1}^{m}({\mathbb{Z}_{2^{\alpha-\lambda}}}^{(2^{\lambda})}\rtimes S 
 _{2^{\lambda}})\bigoplus
\mathcal{H}_{\lambda},& m\geq \lambda +1,\end{cases}$$\par\vspace{.25cm}

\underline{$q\equiv -1~({\rm mod}\,4)$}
$$\mathrm{Aut}(\mathbb{F}_{q}[G])\cong \begin{cases}S_{4}\bigoplus ({\mathbb{Z}_{2}}^{(2^{m}-2)}\rtimes S_{2^{m}-2})
\bigoplus(({\mathrm{SL}_{2}(\mathbb{F}_{q^{2}})\rtimes \mathbb{Z}_{2})}^{(2^{m-2})}
\rtimes S_{2^{m-2}}),&
 m\leq \lambda +1,\\ S_{4}\bigoplus({\mathbb{Z}_{2}}^{(2^{\lambda+1}-2)}\rtimes S_{2^{\lambda+1}-2})
 \bigoplus_{\alpha=\lambda+2}^{m}({\mathbb{Z}_{2^{\alpha-\lambda}}}^{(2^{\lambda})}\rtimes S _{2^{\lambda}})\bigoplus
\mathcal{H}_{\lambda} ,& m\geq \lambda+2,\end{cases}$$ where $\mathcal{H}_{\lambda} = {(\mathrm{SL}_{2}
(\mathbb{F}_{q^{2^{m-\lambda}}})\rtimes \mathbb{Z}_{2^{m-\lambda}})}^{(2^{\lambda-1})}\rtimes S_{2^{\lambda-1}}.$

\end{theorem}\vspace{.25cm}
{\bf{Proof.}}~ We have
$$G =\langle x,y\mid x^{2^{m+1}}=1,y^2=x^{2},y^{-1}xy=x^{2^{m}+1} \rangle.$$ By Lemma \ref{L1}, the non-identity 
normal
subgroups of $G$ are given by
\begin{quote}$(i)$ ~$\langle x^{2^{\alpha}} \rangle,\,\langle x^{2^{\alpha}},x^{2^{\alpha}-1}y
 \rangle,\,0 \leq \alpha \leq m,$\\
$(ii)$~$ \langle x^{2^{\beta}},\,x^{2^{\beta-1}-1}y \rangle,\,
1 \leq \beta \leq m.$
\end{quote}
Also, Lemmas \ref{L5} and \ref{L7} yield that \vspace{.25cm}\\
$\mathcal{S}=
\{(K_{0},\,\langle 1 \rangle,\,K_{1})\} \cup
 \{(K_{1},\,\langle 1 \rangle,\,G/K_{1})\}
\cup \{({K_{2}}^{(\alpha)},\,\langle 1 \rangle,\,G/{K_{2}}^{(\alpha)})\,|\, 0\leq \alpha \leq m\}
\cup \{({K_{3}}^{(\beta)},\,\langle 1 \rangle,\,G/{K_{3}}^{(\beta)})\,|\,1\leq \beta \leq m \}. $
Therefore, $(i)$ follows from Theorem \ref{T11}.\par\vspace{.25cm}
 For each $(N, D/N, A_{N}/N) \in \mathcal{S}$, the corresponding $o(\mathcal{A}_{N},\,D)$ and $|R(\mathcal{A}_{N}/D)|$ 
 in the cases $q\equiv 1~
({\rm mod}\,4)$ or
$q\equiv -1~({\rm mod}\,4)$ are as follows:\pagebreak \par\vspace{.25cm} \noindent
{\bf Case I\,: $q\equiv 1 ~({\rm mod}\,4).$}\vspace{.5cm}\\
\begin{tabular}{|c|c|c|}\hline ($N, \, D/N,\,\mathcal{A}_{N}/N)$& $o(\mathcal{A}_{N},\,D)$&$|R(\mathcal{A}_{N}/D)|$ \\ 
\hline $(K_{1},\,\langle 1 \rangle,\,G/K_{1}),$&1
&1 \\  \hline
$({K_{2}}^{(0)},\,\langle 1 \rangle,\,G/{K_{2}}^{(0)}),$&
$1$ &
$1$\\  \hline
$({K_{2}}^{(\alpha)},\,\langle 1 \rangle,\,G/{K_{2}}^{(\alpha)}),$&
$\begin{cases}2^{\alpha-\lambda},&\alpha\geq \lambda+1,\\1,&\alpha\leq \lambda\end{cases}$ &
$\begin{cases}2^{\lambda-1},&\alpha\geq \lambda+1,\\2^{\alpha-1},&\alpha\leq \lambda\end{cases}$\\ $1\leq \alpha\leq 
m$&&
 \\ \hline
$({K_{3}}^{(\beta)},\,\langle 1 \rangle,\,G/{K_{3}}^{(\beta)}),$& $\begin{cases}
2^{\beta-\lambda},&\beta\geq \lambda+1,\\1,&\beta\leq \lambda
\end{cases}$&
$\begin{cases}2^{\lambda-1},&\beta\geq \lambda+1,\\2^{\beta-1},&\beta\leq \lambda\end{cases}$ \\ $1\leq \beta\leq m$&&
\\ \hline$(K_{0},\, \langle 1 \rangle,\,K_{1}) $&$\begin{cases}2^{m-\lambda},&m\geq \lambda +1,\\1,&m\leq \lambda
\end{cases}$&
$\begin{cases}2^{\lambda-1},&m\geq \lambda +1,\\2^{m-1},&m\leq \lambda\end{cases}$ \\ \hline
\end{tabular}\par\vspace{.5cm}
\par\vspace{.25cm} \noindent
{\bf Case II\,: $q\equiv -1~({\rm mod}\,4).$}\vspace{.25cm}\\

\begin{tabular}{|c|c|c|}\hline ($N, \, D/N,\,\mathcal{A}_{N}/N)$& $o(\mathcal{A}_{N},\,D)$&$|R(\mathcal{A}_{N}/D)|$ \\ 
\hline $(K_{1},\,\langle 1 \rangle,\,G/K_{1}),$&1
&1 \\ \hline
$(K_{2}^{(\alpha)},\,\langle 1 \rangle,\,G/K_{2}^
{(\alpha)}),$& $1$&
$1$ \\  $0\leq \alpha\leq 1$&&\\ \hline
$(K_{2}^{(\alpha)},\,\langle 1 \rangle,\,G/K_{2}^
{(\alpha)}),$& $\begin{cases}2^{\alpha-\lambda},&\alpha\geq \lambda+2,\\2,&\alpha\leq \lambda +1\end{cases}$&
$\begin{cases}2^{\lambda-1},&\alpha\geq \lambda+2,\\2^{\alpha-2},&\alpha\leq \lambda +1\end{cases}$ \\ $2\leq 
\alpha\leq m$&&
\\ \hline $(K_{3}^{(1)},\,
\langle 1 \rangle,\,G/K_{3}^{(1)}),$& $1$
&$1$ \\ \hline $(K_{3}^{(\beta)},\,
\langle 1 \rangle,\,G/K_{3}^{(\beta)}),$& $\begin{cases}2^{\beta-\lambda},&\beta\geq \lambda+2,\\2,&\beta\leq \lambda 
+1
\end{cases}$
&$\begin{cases}2^{\lambda-1},&\beta\geq \lambda+2,\\2^{\beta-2},&\beta\leq \lambda +1\end{cases}$ \\ $2\leq \beta \leq 
m$&& \\ \hline
$(K_{0},\, \langle 1 \rangle,\,K_{1}) $&$\begin{cases}2^{m-\lambda},&m\geq \lambda+2,\\2,&m\leq \lambda 
+1\end{cases}$&
$\begin{cases}2^{\lambda-1},&m\geq \lambda+2,\\2^{m-2},&m\leq \lambda+1\end{cases}$ \\ \hline
\end{tabular}
\vspace{.5cm} \\ Thus Theorem \ref{T2} , with the help of above two tables yield $(ii)$ and $(iii).$~~$\Box$
\bigskip\\
{\bf{Remark:}}~ The above analysis of the structure of $\mathbb{F}_{q}[G],~G$ of type $D_{1},~D_{2},$
provides a method for computing the algebraic structure of $\mathbb{F}_{q}[G],$ for finite group $G$ whose central 
quotient is
Klein four-group. It will thus naturally be of interest to compute the algebraic structure of $\mathbb{F}_{q}[G],~G$ of 
type $D_{i},
~i=3,4,5.$

\end{document}